\documentclass[11pt]{imsart}
\usepackage[margin=1.1in]{geometry}
\RequirePackage[OT1]{fontenc}
\RequirePackage{amsthm,amsmath,natbib}
\RequirePackage[colorlinks,citecolor=blue,urlcolor=blue]{hyperref}
\usepackage{amsfonts}
\usepackage{amssymb}
\usepackage{esint} 
\usepackage[normalem]{ulem}

\newcommand{\real}{\mathbb{R}}
\newcommand{\E}{\text{E}}

\newcommand{\bS}{{S}}
\newcommand{\bV}{{V}}
\newcommand{\bX}{{X}}

\newcommand{\bQ}{{Q}}
\newcommand{\bpi}{{\pi}}
\newcommand{\bdelta}{{\delta}}

\newcommand{\myx}{{\mathcal{X}}}
\newcommand{\myB}{{\mathcal{B}}}

\newtheoremstyle{mytheorem}{}{}{\slshape}{}{\bfseries}{}{2mm}{}
\theoremstyle{mytheorem}
\newtheorem{remark}{Remark}
\newtheorem{theorem}{Theorem}
\newtheorem{corollary}{Corollary}
\newtheorem{proposition}{Proposition}
\newtheorem{definition}{Definition}
\newtheorem{assumption}{Assumption}
\newtheorem{lemma}{Lemma}

\startlocaldefs
\numberwithin{equation}{section}
\theoremstyle{plain}

\theoremstyle{remark}

\endlocaldefs

\begin{document}

\begin{frontmatter}
\title{Central Limit Theorems for Markov Chains from Wasserstein Convergence Rates}
\runtitle{CLTs from Wasserstein rates}

\begin{aug}
	\author{\fnms{Rui} \snm{Jin}}

	
\runauthor{R. Jin}
	
	
\end{aug}

\begin{abstract}
We give sufficient conditions for central limit theorems (CLTs) for additive functionals of Markov chains in terms of quantitative Wasserstein convergence rates, including suitable subgeometric rates. For a given metric $\psi$, we establish CLTs for $\psi$-Lipschitz functions under moment conditions by showing that suitable $1$-Wasserstein convergence rates imply either the Maxwell–Woodroofe projective criterion or convergence of the associated Poisson series. We then extend this framework beyond the $\psi$-Lipschitz setting in two directions. First, by reweighting $\psi$ with a non-negative function $V$, we construct a weighted path metric under which functions with $V$-controlled increments are Lipschitz. We derive convergence bounds in the Wasserstein distance induced by this new metric from corresponding bounds in the Wasserstein distance induced by $\psi$, thereby obtaining CLTs for this broad class of functions. Second, we consider functions admitting integrable increment envelopes and derive CLTs from quantitative $2$-Wasserstein convergence rates. On $\mathbb R^d$, pointwise Sobolev inequalities and polynomial bounds on the gradient provide concrete sufficient conditions for constructing such envelopes. We illustrate the results with nonlinear autoregressive processes and a random walk on the one-dimensional torus exhibiting subgeometric Wasserstein convergence.

\end{abstract}


\end{frontmatter}

\section{Introduction}
Total variation (TV) distance has long been central to the convergence analysis of Markov chains. Classical rates such as geometric and polynomial ergodicity are typically formulated in TV, and there is a mature toolkit for obtaining qualitative and quantitative TV bounds; see, for example, \citet{rose:1995,jarner:robert:2002}. Such rates are closely connected to mixing properties and therefore provide a natural route to central limit theorems (CLTs) for additive functionals of Markov chains; see \citet{jones:2004} and references therein. For instance, geometric ergodicity implies exponentially fast strong mixing under standard conditions \citep{chan:geyer:1994}, allowing classical mixing-based CLTs to be applied.

However, convergence analyses in total variation, particularly those based on classical single-step drift-and-minorization arguments, can become difficult to implement and may produce conservative quantitative bounds in complex or high-dimensional models. This limitation has been documented in convergence complexity studies; see, for example,  \citet{qin:hobe:2021aap}. Wasserstein distances provide an alternative that is often more amenable to direct coupling arguments and to problem-adapted choices of state-space geometry. A substantial literature has developed geometric and subgeometric Wasserstein convergence bounds using contraction, Lyapunov conditions, weighted distances, and related coupling techniques; see, among others, \citet{hairer:mattingly:scheutzow:2011,butkovsky:2014,durmus:moulines:2015,eberle:majka:2019,qin:hobe:2019}. Thus, quantitative Wasserstein convergence bounds are now available for a broad range of Markov chains, including settings in which a comparable TV analysis may be substantially less tractable.

Despite this increasingly flexible toolkit for quantitative Wasserstein convergence, comparatively less work has used such rates directly to establish CLTs. 
An early result of \citet{walczuk:2008} established a CLT for additive functionals of Markov processes under the assumption that the process is asymptotically contractive in the Wasserstein metric. \citet{komorowski:walczuk:2012} subsequently proved CLTs for additive functionals of continuous-time Markov processes under strong contraction in an appropriate Wasserstein metric. \citet{gulgowski:hille:szarek:ziemlanska:2019} obtained a CLT for non-stationary Markov chains under a spectral-gap condition. For random walks on compact groups, \citet{Borda2021} showed that summability of suitable Wasserstein distances to Haar measure yields CLTs, as well as stronger limit theorems, for H\"older observables. More recently, \citet{czapla:horbacz:wojewodka:2024} considered continuous-time Markov--Feller processes under exponential ergodicity in the bounded-Lipschitz metric together with a continuous Foster--Lyapunov condition. 

In summary, quantitative Wasserstein convergence bounds are now available for a wide range of Markov chains and can often be more robust to increasing dimension, particularly when they are obtained through contractive couplings and problem-adapted cost functions. However, existing Wasserstein-based CLTs generally impose regularity conditions on the observable with respect to the underlying state-space metric, most commonly Lipschitz or H\"older continuity.
This can be restrictive in applications when functions of interest do not satisfy the required regularity for the metric that induces the Wasserstein distance in which convergence is available. The resulting question is how to use an available Wasserstein convergence rate to obtain CLTs for functions outside these regularity classes. Extending the Wasserstein-based approach to such functions requires additional geometric or analytic structure.

This paper develops a framework for proving central limit theorems for additive functionals of discrete-time Markov chains directly from quantitative Wasserstein convergence bounds. Let $(r(n))_{n\geq 1}$ denote a quantitative Wasserstein convergence rate sequence. Our main results are as follows.
\begin{itemize}
\item \textbf{Part I: CLTs for Lipschitz functions from $1$-Wasserstein convergence.} 
We show that, under suitable moment conditions, the summability condition
\[
\sum_{n\geq 1}\frac{r(n)}{\sqrt n}<\infty
\]
is sufficient to establish CLTs for square-integrable Lipschitz functions. The proof verifies  the Maxwell--Woodroofe projective criterion directly from the Wasserstein rate, without appealing to total variation mixing, spectral gaps, or irreducibility. Consequently, the results also apply in settings where standard TV-based arguments are unavailable or difficult to use, including chains that do not converge in total variation and chains with subgeometric Wasserstein convergence.

\item \textbf{Part II: Extending CLTs beyond Lipschitz functions.}
We then develop two approaches for obtaining CLTs for functions that need not be Lipschitz with respect to the original metric $\psi$.
\begin{enumerate}
\item[(a)] \emph{Weighted path-metric extension.}
We reweight the state-space geometry using a non-negative function $V$. Starting from $\psi$ and $V$, we introduce a weighted cost and construct the corresponding path metric. Using a weighted H\"older coupling argument, related to that in~\cite{kulik:2018}, we then derive convergence bounds in the Wasserstein distance induced by this path metric from the original Wasserstein bounds, with the rate $r(n)$ transformed to $r(n)^\alpha$ for $\alpha\in(0,1)$. Part~I can then be applied with the new metric, yielding CLTs for functions that are Lipschitz with respect to the weighted path metric and whose increments may grow with $V$.
\item[(b)] \emph{Increment-envelope extension.} For a function $g$ that need not be $\psi$-Lipschitz, we instead assume that its increments $|g(x)-g(y)|$ admit a suitable integrable envelope. Combining this increment envelope with a quantitative $2$-Wasserstein convergence rate and appropriate integrability conditions, we verify the Maxwell--Woodroofe criterion and hence obtain a CLT for $g$ when
$\sum_{n\geq 1}{r(n)}/{\sqrt n}<\infty$.

\end{enumerate}
\end{itemize}

Section~\ref{preliminary} provides the notation and definitions used in the paper. Section~\ref{main_result} presents CLTs from quantitative Wasserstein rates and their extensions. 
Section~\ref{application} provides two examples illustrating the application of these results. For readability, most proofs and ancillary technical material are deferred to the Appendix.

\section{Preliminaries}\label{preliminary}

\subsection{General state space Markov chains}

Let $\myx$ be a Polish space with Borel $\sigma$-field $\myB$, and let
$\bQ$ be a Markov transition kernel on $(\myx,\myB)$. Thus,
$\bQ(x,\cdot)$ is a probability measure for every $x\in\myx$, and
$x\mapsto \bQ(x,C)$ is $\myB$-measurable for every $C\in\myB$.
Let $ (\bX_n)_{n\geq 0} $
be a time-homogeneous Markov chain with transition kernel $\bQ$.

For a probability measure $\mu$ on $(\myx,\myB)$, define
\[
    (\mu\bQ)(C)
    :=
    \int_{\myx}\bQ(x,C)\,\mu(dx),
    \qquad C\in\myB.
\]
Write $\bQ^n$ for the $n$-step transition kernel, with $\bQ^0$
denoting the identity kernel. Then $\bdelta_x\bQ^n=\bQ^n(x,\cdot)$
is the distribution of $\bX_n$ when $\bX_0=x$. We write
$\mathbb P_\mu$ and $\E_\mu$ for probability and expectation of the chain when
$\bX_0\sim\mu$, and abbreviate $\mathbb P_x:=\mathbb P_{\bdelta_x},
    \E_x:=\E_{\bdelta_x}$.

Assume that $\bQ$ admits an invariant probability measure $\bpi$, that is, 
  $  \bpi(C) = \int_{\myx}\bQ(x,C)\,\bpi(dx), C\in\myB$.
When $\bX_0\sim\bpi$, the chain is stationary.

For a probability measure $\mu$ and a $\mu$-integrable function $g$,
write  $\mu(g):=\int_{\myx}g(x)\,\mu(dx)$.
Then, let
\[
    L^2(\bpi)
    :=
    \left\{
        g:\myx\to\real:
        g\text{ is $\myB$-measurable and }\bpi(g^2)<\infty
    \right\},
\]
and
\[
    L_0^2(\bpi)
    :=
    \left\{
        g\in L^2(\bpi):\bpi(g)=0
    \right\}.
\]
For $g\in L^2(\bpi)$, write
    $\|g\|
    :=
    \left(
        \int_{\myx}g^2(x)\,\bpi(dx)
    \right)^{1/2}$.
Unless otherwise specified, $\|\cdot\|$ denotes the $L^2(\bpi)$
norm.

The same symbol $\bQ$ denotes the associated Markov operator  on
$L^2(\bpi)$. Define
\begin{equation*}
\bQ g(x) = \int_{\myx} g(y) \bQ(x; dy) \;\;\text{for any $g \in L^2(\bpi)$.}
\end{equation*}
It can be shown that $\bQ$ is a contraction, in the sense that  $\| \bQ g \| \leq \| g \|$ for any $g \in L^2(\bpi)$. 

Throughout this paper, when the chain is said to be ergodic with respect to $\pi$, we mean that the stationary chain under $\mathbb P_\pi$ is ergodic.

\subsection{Wasserstein distances and Lipschitz observables}\label{sec:wasserstein-prelim}

Let $(\myx, \psi)$ be a Polish metric space.  Further, let ${\mathcal{P}}(\myx)$ be the set of probability measures on $(\myx, \myB)$ and $\bdelta_x$ the point mass at  $x$. 
For probability measures $\mu$ and $\nu$, let $\mathcal C(\mu,\nu)$ denote
the set of probability measures on $\mathcal X\times\mathcal X$ with
marginals $\mu$ and $\nu$.  Its elements are couplings of $\mu$ and $\nu$.
For $p\geq1$, define the $p$--Wasserstein distance (induced by $\psi$)
\begin{equation}\label{eq:Wp}
  W_{\psi,p}(\mu,\nu)
  :=\left\{
  \inf_{\Gamma\in\mathcal C(\mu,\nu)}
  \int_{\mathcal X\times\mathcal X}\psi(x,y)^p\,\Gamma(dx,dy)
  \right\}^{1/p}.
\end{equation}
The results from Section~\ref{suff_mapx} use the case $p=1$, for which we write
$
  W_\psi:=W_{\psi,1}.
$
The later increment-envelope argument uses $p=2$.
To specify the natural domain of these distances, let
\begin{equation}\label{w_space}
  \mathcal P_\psi^p
  :=\left\{\mu:
  \int\psi(x_0,x)^p\,\mu(dx)<\infty
  \text{ for some }x_0\in\mathcal X\right\}.
\end{equation}
When $\mathcal X=\mathbb R^d$ and $\psi(x,y)=|x-y|$, where $|\cdot|$ is the
Euclidean norm, we use the standard abbreviations $W_p:=W_{\psi,p}$ and
$\mathcal P_p:=\mathcal P_\psi^p$.

The dual class associated with $W_\psi$ consists of Lipschitz functions.  For
a real-valued function $g$, define
\[
  \operatorname{Lip}_\psi(g)
  :=\sup_{x\neq y}\frac{|g(x)-g(y)|}{\psi(x,y)},
\]
and let
\begin{equation}\label{eq:Gpsi}
  \mathcal G_\psi:=\{g:\operatorname{Lip}_\psi(g)\leq1\}.
\end{equation}
For $\mu,\nu\in\mathcal P_\psi^1$, Kantorovich--Rubinstein duality gives
\begin{equation}\label{dual_wass}
  W_\psi(\mu,\nu)
= \sup_{g \in \mathcal{G}_{\psi}} \left|\int_{\myx} g(x) \mu(dx) - \int_{\myx} g(x) \nu(dx)  \right|   =\sup_{g\in\mathcal G_\psi}|\mu(g)-\nu(g)|.
\end{equation}

\subsection{Martingale approximation and martingale CLTs}\label{mapprox_mclt}
For $g\in L_0^2(\bpi)$ and $n\geq1$, define the additive functional
\[
    \bS_n(g)
    :=
    \sum_{k=0}^{n-1}g(\bX_k),
\]
and let
\[
    {\mathcal F}_n
    :=
    \sigma(\bX_0,\ldots,\bX_n),
    \qquad n\geq0.
\]
Unless otherwise stated, all expectations and martingale approximations
in this subsection are taken under the stationary law
$\mathbb P_{\bpi}$. Our approach is to approximate $\bS_n(g)$ by a martingale, with an
approximation error that is negligible on the $\sqrt n$ scale.

\begin{definition}\label{mart_approx}
We say that $\bS_n(g)$ admits an $L^2$ martingale approximation under
$\mathbb P_{\bpi}$ if there exists a square-integrable martingale
\[
    \bigl( M_n,{\mathcal F}_n\bigr)_{n\geq0},
    \qquad  M_0=0,
\]
whose increment sequence
\[
    \bigl( M_k- M_{k-1}\bigr)_{k\geq1}
\]
is stationary under $\mathbb P_{\bpi}$, and such that
\[
    \E_{\bpi}
    \left[
        \bigl(\bS_n(g)- M_n\bigr)^2
    \right]
    =
    o(n)
    \qquad\text{as }n\to\infty.
\]
\end{definition}
Let     $R_n := \bS_n(g)- M_n$, the $L^2$ martingale approximation implies
\begin{equation}\label{martingale_approx}
    \frac{\bS_n(g)}{\sqrt n}
    =
    \frac{ M_n}{\sqrt n}
    +
    \frac{ R_n}{\sqrt n}
    =
    \frac{ M_n}{\sqrt n}
    +
    o_{\mathbb P_{\bpi}}(1).
\end{equation}

We use two related sufficient conditions for the existence of such an
approximation. The first constructs the martingale explicitly through
the Poisson equation. The second is the more flexible
Maxwell--Woodroofe projective condition, which controls the cumulative
conditional means of the partial sums. First, consider the Poisson equation
\begin{equation}\label{poiseq}
    h-\bQ h=g.
\end{equation}
If \eqref{poiseq} has a solution $h\in L^2(\bpi)$, then
\begin{equation*}\label{eq:poisson-decomposition}
\begin{split}
    \bS_n(g)
    &=
    \sum_{k=1}^{n}
    \left\{
        h(\bX_k)-\bQ h(\bX_{k-1})
    \right\}
    +
    h(\bX_0)-h(\bX_n) \\
    &=
    M_n+ R_n.
\end{split}
\end{equation*}

A convenient sufficient condition for the existence of a solution to
\eqref{poiseq} is the Poisson-series condition
\begin{equation}\label{solvePE}
    \sum_{k=0}^{\infty}\|\bQ^k g\|<\infty.
\end{equation}
Under \eqref{solvePE}, the series $ h
    :=
    \sum_{k=0}^{\infty}\bQ^k g$
converges in $L^2(\bpi)$. Since $\bQ$ is a contraction on
$L^2(\bpi)$, it follows that $h-\bQ h=g$.
The Poisson-series condition is convenient because it constructs an
explicit solution to the Poisson equation. It can, however, be stronger
than necessary for the existence of a martingale approximation. A second
approach, due to \citet{maxwell:woodroofe:2000}, instead controls the
partial sums of the iterates $\bQ^k g$.
For $n\geq 1$, define
\[
    \bV_n g
    :=
    \sum_{k=0}^{n-1}\bQ^k g.
\]

 The following
projective condition gives a martingale approximation without requiring
an explicit solution to the Poisson equation.

\begin{theorem}[Maxwell--Woodroofe]
\label{sufficient_martingale_approx}
Suppose that the chain is ergodic with respect to $\bpi$, and let
$g\in L_0^2(\bpi)$. If
\begin{equation}\label{suff_cond_clt}
    \sum_{n=1}^{\infty}
    n^{-3/2}\|\bV_n g\|
    <
    \infty,
\end{equation}
then $\bS_n(g)$ admits an $L^2$ martingale approximation under
$\mathbb P_{\bpi}$.
\end{theorem}

The Maxwell--Woodroofe condition is potentially weaker than~\eqref{solvePE}: it
controls the growth of the partial sums $\bV_n g$ and does not require
absolute convergence of the Poisson series
$\sum_{k\geq0}\bQ^k g$ in $L^2(\bpi)$. 
The preceding results concern martingale approximation under the
stationary law. Our main results also require a CLT when the chain is
started from a fixed state. For this purpose, we use the following
Markov chain consequence of the quenched weak invariance principle of
\citet{cuny:merlevede:2014}.

\begin{theorem}[Cuny--Merlev\`ede]
\label{clt_anypoint}
Suppose that the chain is ergodic with respect to $\bpi$, and let
$g\in L_0^2(\bpi)$. If \eqref{suff_cond_clt} holds, then
\[
    \sigma^2(g)
    :=
    \lim_{n\to\infty}
    \frac{1}{n}
    \E_{\bpi}
    \left[
        \bS_n(g)^2
    \right]
\]
exists and is finite. Moreover, for $\bpi$-almost every
$x\in\myx$,
\begin{equation}\label{eq:quenched-clt}
    \frac{\bS_n(g)}{\sqrt n}
    \Rightarrow
    N\!\left(0,\sigma^2(g)\right)
    \qquad\text{under }\mathbb P_x
\end{equation}
as $n\to\infty$.
\end{theorem}

\section{Main results}\label{main_result}

\subsection{CLT for $\psi$--Lipschitz functions}\label{suff_mapx}

 Let
$\Lambda:\mathcal X\to[0,\infty)$ be measurable, and suppose that
\begin{equation}\label{w_subge}
  W_\psi(\delta_xQ^n,\pi)
  \leq \Lambda(x)r(n),
  \qquad x\in\mathcal X,\quad n\geq1.
\end{equation}
For $g\in\mathcal G_\psi\cap L_0^2(\pi)$, we then have
\begin{equation}\label{eq:key-L2-decay}
  |Q^ng(x)|
  =|\delta_xQ^n(g)-\pi(g)|
  \leq \Lambda(x)r(n).
\end{equation}
Consequently,
\begin{equation}\label{eq:key-L2-decay-norm}
  \|Q^ng\|\leq\|\Lambda\|r(n).
\end{equation}
This inequality reveals the distinct roles of the assumptions:  the rate $r(n)$ determines whether the projective sums are
summable, while the moment condition on $\Lambda$ converts a pointwise bound
into an $L^2(\pi)$ bound.

\begin{theorem}\label{main_thm}
Suppose~\eqref{w_subge} and $\Lambda\in L^2(\pi)$ hold, and assume
\begin{equation}\label{rate_fun}
  \sum_{n=1}^{\infty}\frac{r(n)}{\sqrt n}<\infty.
\end{equation}
Then every $g\in\mathcal G_\psi\cap L_0^2(\pi)$ satisfies
\begin{equation*}
  \sum_{n=1}^{\infty}n^{-3/2}
  \left\|\sum_{k=0}^{n-1}Q^kg\right\|<\infty.
\end{equation*}
If the stationary chain is ergodic, then, for $\pi$-almost every $x$,
\[
  \frac{S_n(g)}{\sqrt n}
  \Rightarrow N(0,\sigma^2(g))
  \qquad\text{under }\mathbb P_x.
\]

\end{theorem}
Condition~\eqref{rate_fun} is readily verified for a broad range of commonly used convergence rates. In particular, it holds for the polynomial rate  $r(n)=n^{-\beta}$ whenever $\beta>1/2$, and for the rate $r(n)=\rho^{n^\gamma}$ with
$0<\rho<1$ and $\gamma>0$. If the rate is summable without the factor $n^{-1/2}$, this set of conditions gives
a stronger Poisson-series conclusion. 

\begin{corollary}\label{check:solvePE}
Suppose~\eqref{w_subge} and $\Lambda\in L^2(\pi)$ hold.  If
\begin{equation}\label{eq:poisson-rate}
  \sum_{n=1}^{\infty}r(n)<\infty,
\end{equation}
then every $g\in\mathcal G_\psi\cap L_0^2(\pi)$ satisfies~\eqref{solvePE}.
\end{corollary}

We next illustrate how  
assumptions~\eqref{rate_fun} and~\eqref{eq:poisson-rate} relate to the one-step contraction. 
Define the Wasserstein
contraction coefficient
\[
  \Delta_\psi(Q)
  :=\sup_{x\neq y}
  \frac{W_\psi(\delta_xQ,\delta_yQ)}{\psi(x,y)}.
\]

\begin{proposition}\label{MA_GC}
Suppose $\Delta_\psi(Q)<1$ and $\pi\in\mathcal P_\psi^2$.  Then
\eqref{w_subge}--\eqref{rate_fun} hold with a geometric rate, and
\eqref{solvePE} holds for every
$g\in\mathcal G_\psi\cap L_0^2(\pi)$.  If the stationary chain is ergodic,
the quenched CLT in Theorem~\ref{main_thm} also holds.
\end{proposition}

\subsection{Weighted path-metric geometry and convergence rates}\label{weight_lip}

{
Theorem~\ref{main_thm} from Section~\ref{suff_mapx} only applies to $\psi$-Lipschitz functions.  This excludes, for example, unbounded
functions when $\psi$ is bounded.  Even for an unbounded metric $\psi$,
common functions need not be $\psi$-Lipschitz.   For illustration,  let $g(x) = x^2$ on $\real$ with the Euclidean metric, we then have
  $$|g(x)-g(y)|=|x-y|\,|x+y|.$$
The factor $|x+y|$ is unbounded, so $g$ is not Euclidean
Lipschitz.  At the same time, we observe that 
\[
|g(x)-g(y)|
\leq |x-y|(|x|+|y|) \,.
\]
This suggests that if  the distance between $x$ and $y$ is enlarged from $|x-y|$ to a distance that is at least  $|x-y|(|x|+|y|)$,
then $g(x)=x^2$ becomes Lipschitz with respect to the enlarged distance.
Motivated by this observation, we develop a weighted path-metric geometry and establish CLTs for functions that become Lipschitz under the resulting metric.

}

{
\subsubsection{Weighted costs and their path metrics}\label{weighted_cost}
Let $V:\mathcal X\to[0,\infty)$ be a continuous weight function.  Fix $\alpha\in(0,1)$ and $A>0$,  define the weighted cost 
\begin{equation}\label{eq:weighted-cost}
  c_{\alpha,A}(x,y)
  :=\psi(x,y)^\alpha
  \{A+V(x)+V(y)\}^{1-\alpha}.
\end{equation}

We use this cost to specify state-dependent control on the increment $|g(x)-g(y)|$.
Accordingly, define the weighted   H\"older class by
\begin{equation}\label{eq:weighted-class}
  \mathcal G_{\psi,V}^{(\alpha,A)}
  :=\left\{g:
  |g(x)-g(y)|\leq c_{\alpha,A}(x,y)
  \text{ for all }x,y\in\mathcal X\right\}.
\end{equation}

To use this increment control together with Wasserstein convergence of the chain, we seek a metric whose Lipschitz class coincides with $\mathcal G_{\psi,V}^{(\alpha,A)}$. The cost $c_{\alpha,A}$ itself need not be a metric because it may fail the triangle inequality. We therefore construct the path metric
\begin{equation}\label{eq:simple-path-metric}
  d_{\alpha,A}(x,y)
  :=\inf_{m\geq1}\ \inf_{x=z_0,\ldots,z_m=y}
  \sum_{i=1}^m c_{\alpha,A}(z_{i-1},z_i),
\end{equation}
where the inner infimum is over all finite paths from $x$ to $y$.  
The following proposition shows that $ d_{\alpha,A}$ is a Polish metric on $\myx$, and that  replacing $c_{\alpha,A}$ by $d_{\alpha,A}$ in~\eqref{eq:weighted-class} does not change the function class $\mathcal G_{\psi,V}^{(\alpha,A)}$.
\begin{proposition}
\label{lem:path-metric-geometry}
For every $x,y\in\mathcal X$,
\begin{equation}\label{eq:path-lower-upper}
  A^{1-\alpha}\psi(x,y)^\alpha
  \leq d_{\alpha,A}(x,y)
  \leq c_{\alpha,A}(x,y).
\end{equation}
Consequently, $d_{\alpha,A}$ is a Polish metric that generates the same
topology and Borel $\sigma$-field as $\psi$.  

Moreover, the weighted increment condition in~\eqref{eq:weighted-class}
is equivalent to being  $1$-Lipschitz with respect to $d_{\alpha,A}$, that is 
\begin{equation}\label{eq:weighted-class-equivalence}
  \mathcal G_{\psi,V}^{(\alpha,A)}
  =\{g:\operatorname{Lip}_{d_{\alpha,A}}(g)\leq1\}.
\end{equation}
\end{proposition}

Thus, $\mathcal G_{\psi,V}^{(\alpha,A)}$ is exactly the Lipschitz class under $d_{\alpha,A}$. It remains to control convergence of the chain in the corresponding Wasserstein distance $W_{d_{\alpha,A}}$. Since the available convergence bound is expressed in terms of $W_\psi$, we next relate $W_{d_{\alpha,A}}$ directly to $W_\psi$. The following theorem shows that finite first $V$-moments are sufficient for this transfer. The underlying
weighted H\"older coupling argument appears in
\cite{kulik:2018}; the path-metric representation above
yields the following Wasserstein distance form.

\begin{theorem}\label{thm:direct-lift}
If $\mu(V) < \infty$ and $ \nu(V)<\infty$, then
\begin{equation}\label{eq:direct-lift}
  W_{d_{\alpha,A}}(\mu,\nu)
  \leq
  \{A+\mu(V)+\nu(V)\}^{1-\alpha}
  W_\psi(\mu,\nu)^\alpha.
\end{equation}
\end{theorem}

Theorem~\ref{thm:direct-lift} itself requires only first moments of $V$. To apply it with $\mu=\delta_xQ^n$ and $\nu=\pi$, however, we need to control $Q^nV(x)$ along the chain. Moreover, the resulting starting-point envelope must belong to $L^2(\pi)$ in order to apply Theorem~\ref{main_thm}. The following proposition gives sufficient conditions for both requirements.

\begin{proposition}\label{prop:direct-lift}
Suppose~\eqref{w_subge} holds with $\Lambda\in L^2(\pi)$.  Assume that
$V\in L^2(\pi)$ and, for some $C_V<\infty$,
\begin{equation}\label{eq:direct-moment-iterate}
  \sup_{n\geq1}Q^nV(x)\leq C_V\{1+V(x)\},
  \qquad x\in\mathcal X.
\end{equation}
Then there is a finite constant $C$ such that
\begin{equation}\label{eq:direct-transferred-envelope}
  W_{d_{\alpha,A}}(\delta_xQ^n,\pi)
  \leq \widetilde\Lambda_\alpha(x)r(n)^\alpha,
  \qquad
  \widetilde\Lambda_\alpha(x)
  :=C\Lambda(x)^\alpha\{1+V(x)\}^{1-\alpha},
\end{equation}
and $\widetilde\Lambda_\alpha\in L^2(\pi)$.

If, in addition,
\begin{equation}\label{eq:weighted-rate-summability}
  \sum_{n=1}^{\infty}\frac{r(n)^\alpha}{\sqrt n}<\infty,
\end{equation}
then every
$g\in\mathcal G_{\psi,V}^{(\alpha,A)}\cap L_0^2(\pi)$ satisfies the
Maxwell--Woodroofe condition.  If the stationary chain is ergodic, the
quenched CLT~\eqref{eq:quenched-clt} holds for every such $g$ and for $\pi$-almost every initial state.
\end{proposition}

Condition~\eqref{eq:direct-moment-iterate} is immediate if $V$ satisfies a geometric drift inequality. Indeed, if \[ QV\leq\lambda V+b_V, \qquad 0\leq\lambda<1,\quad b_V<\infty, \] then \[ Q^nV(x) \leq\lambda^nV(x)+\frac{b_V}{1-\lambda} \leq V(x)+\frac{b_V}{1-\lambda}. \]

\subsubsection{Weighted H\"older classes with separate exponents}\label{holder_class}

The above weighted H\"older class $ \mathcal G_{\psi,V}^{(\alpha,A)}$ fixes the exponent on the weight function $V$ at $1-\alpha$. It arises from  H\"older's inequality in the proof of Theorem~\ref{thm:direct-lift}. When higher moments of $V$ are available, the exponent on $V$ can instead be chosen separately from $\alpha$, in particular allowing $\gamma > 1 - \alpha$, while the transferred rate remains $r(n)^\alpha$.

 For $\alpha \in (0, 1)$ and  $\gamma >0$,  define
\begin{align}
  c_{\alpha,\gamma,A}(x,y)
  &:=\psi(x,y)^\alpha
  \{A+V(x)+V(y)\}^{\gamma},
  \label{eq:weighted-holder-cost}\\
  d_{\alpha,\gamma,A}(x,y)
  &:=\inf_{m\geq1}\ \inf_{x=z_0,\ldots,z_m=y}
  \sum_{i=1}^m c_{\alpha,\gamma,A}(z_{i-1},z_i),
  \label{eq:path-metric} \\
  \mathcal G_{\psi,V}^{(\alpha,\gamma,A)}
  &:=\left\{g:
  |g(x)-g(y)|\leq c_{\alpha,\gamma,A}(x,y)
  \text{ for all }x,y\in\mathcal X\right\}.
  \label{def:weighted-holder}
\end{align}

Let $s= \gamma / (1 - \alpha) $. We next extend results of Proposition~\ref{lem:path-metric-geometry} and Theorem~\ref{thm:direct-lift} to this
 general cost.

\begin{theorem}\label{thm:weighted-holder-lift}
For every $x,y\in\mathcal X$,
\begin{equation}\label{eq:weighted-holder-path-bounds}
  A^\gamma\psi(x,y)^\alpha
  \leq d_{\alpha,\gamma,A}(x,y)
  \leq c_{\alpha,\gamma,A}(x,y).
\end{equation}
Consequently, $d_{\alpha,\gamma,A}$ is a Polish metric with the same topology
and Borel sets as $\psi$, and
\[
  \mathcal G_{\psi,V}^{(\alpha,\gamma,A)}
  =\{g:\operatorname{Lip}_{d_{\alpha,\gamma,A}}(g)\leq1\}.
\]
If $\mu(V^s)  <\infty$ and $\nu(V^s)<\infty$, then
\begin{equation}\label{eq:wh-transfer}
  W_{d_{\alpha,\gamma,A}}(\mu,\nu)
  \leq C_{\alpha,\gamma}
  W_\psi(\mu,\nu)^\alpha
  \{A^s+\mu(V^s)+\nu(V^s)\}^{1-\alpha},
\end{equation}
where
\[
  C_{\alpha,\gamma}=c_s^{1-\alpha},
  \qquad
  c_s=
  \begin{cases}
    1,&0<s\leq1,\\
    3^{s-1},&s>1.
  \end{cases}
\]
\end{theorem}

Under the corresponding $V^s$-moment control along the chain, the CLT result takes the same form as Proposition~\ref{prop:direct-lift}.

\begin{proposition}
\label{prop:weighted-holder-lift}
Suppose~\eqref{w_subge} holds with $\Lambda\in L^2(\pi)$. Assume that
$V^s\in L^2(\pi)$ and, for some finite constant $C_{V,s}$,
\begin{equation}\label{eq:weighted-moment-iterate}
\sup_{n\geq1}Q^nV^s(x)
\leq C_{V,s} \{1+V(x)^s\},
\qquad x\in\mathcal X.
\end{equation}

Then there is a finite constant $C$ such that
\begin{equation}\label{eq:direct-transferred-envelope-general}
W_{d_{\alpha,\gamma,A}}(\delta_xQ^n,\pi)
\leq \widetilde\Lambda_{\alpha,\gamma}(x)r(n)^\alpha,
  \qquad
  \widetilde\Lambda_{\alpha,\gamma}(x)
  :=C\Lambda(x)^\alpha\{1+V(x)^s\}^{1-\alpha},
\end{equation}
and $\widetilde\Lambda_{\alpha,\gamma}\in L^2(\pi)$.

If, in addition,~\eqref{eq:weighted-rate-summability} holds, then every
$g\in\mathcal G_{\psi,V}^{(\alpha,\gamma,A)}\cap L_0^2(\pi)$ satisfies
the Maxwell--Woodroofe condition. If the stationary chain is ergodic,
the quenched CLT~\eqref{eq:quenched-clt} holds for every such $g$ and for $\pi$-almost every initial state.
\end{proposition}

}

\subsection{Function-specific increment envelopes under
\texorpdfstring{$W_{\psi,2}$}{W2}}\label{sobolev_class}

The geometric approach in Section~\ref{weight_lip} is appropriate when the increments of a class of  functions can be controlled using a common weight function $V$. For some functions, however, the  increment bound depends on the function itself rather than through a common weight. A familiar Euclidean example illustrates this distinction. For
$g\in C^1(\mathbb R^d)$,  the fundamental theorem of calculus gives
\begin{equation*}\label{eq:line-segment-motivation}
  |g(x)-g(y)|
  =|\int_0^1\nabla g\{y+t(x-y)\}^{T}(x-y)\,dt| \leq |x - y| \int_0^1 |\nabla g\{y+t(x-y)\} | \,dt.
\end{equation*}
Thus, the increment of $g$ is controlled by the Euclidean distance $|x-y|$
multiplied by a quantity determined by $g$ itself.
Motivated by this bound, we seek a function-specific envelope $G_g$ that controls this quantity, for example,
\[
\int_0^1 |\nabla g\{y+t(x-y)\} | \,dt \leq G_g(x) + G_g(y) \,.
\]
Such a bound would give 
\[
|g(x)-g(y)|
\leq |x-y|\{G_g(x)+G_g(y)\}.
\]
This suggests an alternative extension of results from Section~\ref{suff_mapx}. Rather than re-weighting the geometry of the state space to expand the range of Lipschitz functions, we keep the metric $\psi$ 
and study CLTs for functions admitting integrable increment envelopes. We formalize this class of functions as follows.

\begin{definition}\label{def:increment-envelope}
Let $q\geq2$.  A function $g\in L_0^2(\bpi)$ is said to admit an
$L^q(\bpi)$ \emph{increment envelope} if there exist a measurable function
$G_g:\mathcal X\to[0,\infty]$ and a constant $C_g<\infty$ such that
\begin{equation}\label{eq:increment-envelope}
  |g(x)-g(y)|
  \leq C_g\psi(x,y)\{G_g(x)+G_g(y)\},
  \qquad x,y\in\mathcal X,
\end{equation}
and $G_g\in L^q(\bpi)$.
\end{definition}
For functions studied in Sections~\ref{suff_mapx} and~\ref{weight_lip}, their increment bound involves only the underlying metric $\psi$ or the corresponding path metric $d_{\alpha, A}$. In contrast,   for a function with an $L^q(\bpi)$ {increment envelope}, \eqref{eq:increment-envelope} contains the additional factor $G_g(x)+G_g(y)$.  Thus, for an optimal coupling $\gamma$ of
$\delta_xQ^n$ and $\pi$, the relevant quantity becomes
\[
\int
\psi(y,z)\{G_g(y)+G_g(z)\}\gamma(dy,dz).
\]
By the Cauchy--Schwarz inequality,
\[
\int
\psi(y,z)\{G_g(y)+G_g(z)\} \gamma(dy,dz)
\leq
\left\{\int\psi(y,z)^2\gamma(dy,dz)\right\}^{1/2}
\left\{\int\{G_g(y)+G_g(z)\}^2\gamma(dy,dz)\right\}^{1/2}.
\]
The first factor is controlled by the  $2$-Wasserstein distance $W_{\psi, 2}$, while the second is controlled through moments of the increment envelope.
This leads to the following $W_{\psi,2}$ convergence condition on the
Markov chain.

\begin{assumption}\label{ass:Wpsi2-envelope}
The invariant law satisfies $\pi\in\mathcal P_\psi^2$.  There are a rate
sequence $r$ and a measurable function
$\Lambda:\mathcal X\to[0,\infty)$ such that
\begin{equation}\label{eq:Wpsi2-envelope-rate}
  W_{\psi,2}(\delta_xQ^n,\pi)
  \leq \Lambda(x)r(n),
  \qquad x\in\mathcal X,\quad n\geq1,
\end{equation}
with
\[
  \sum_{n=1}^{\infty}\frac{r(n)}{\sqrt n}<\infty.
\]
\end{assumption}

The next theorem establishes CLTs for functions with  increment envelopes  based on the  $W_{\psi,2}$ convergence bound.  
\begin{theorem}\label{thm:increment-envelope-lift}
Suppose Assumption~\ref{ass:Wpsi2-envelope} holds.  Let
$a,b\in(1,\infty)$ satisfy $a^{-1}+b^{-1}=1$.  Assume that
$\Lambda\in L^{2a}(\bpi)$ and that $g$ admits an
$L^{2b}(\bpi)$ increment envelope.  Then, for every $k\geq1$,
\begin{equation}\label{eq:L2rate-envelope}
  \|\bQ^kg\|_2
  \leq 2C_g\|\Lambda\|_{2a}\|G_g\|_{2b} r(k).
\end{equation}
Consequently, $g$ satisfies the Maxwell--Woodroofe condition.  If the
stationary chain is ergodic, then the quenched CLT
\eqref{eq:quenched-clt} holds for $g$ and for $\bpi$-almost every initial
state.
\end{theorem}

\begin{remark}\label{rem:symmetric-L4}
 Taking $a=b=2$ gives the particularly simple sufficient
conditions
\[
  \Lambda,G_g\in L^4(\bpi),
\]
and the bound
\[
  \|\bQ^kg\|_2
  \leq2C_g\|\Lambda\|_4\|G_g\|_4r(k).
\]
 We use the symmetric case below in Sections~\ref{sobolev_sufficient} and~\ref{poly_sufficient} since it gives  simple conditions to verify in applications.
\end{remark}

We now return to the Euclidean setting that motivated the increment-envelope condition. Suppose $\mathcal X=\mathbb R^d$ with
$\psi(x,y)=|x-y|$, where $|\cdot|$ is the Euclidean norm.  Then,
$W_{\psi,2}=W_2$.    The next two sections give concrete sufficient conditions, based on derivative information, for constructing the increment envelopes required by Definition~\ref{def:increment-envelope}.

\subsubsection{A Sobolev construction}\label{sobolev_sufficient}

The motivating example of Section~\ref{sobolev_class} involves the
gradient at every point between $x$ and $y$.  To obtain an  inequality
of the form~\eqref{eq:increment-envelope}, we use a pointwise Sobolev inequality that bounds this increment in terms of the Hardy--Littlewood maximal function evaluated at $x$ and $y$.

For a locally integrable function $h:\mathbb R^d\to\mathbb R$, let
$B(x,t)$ be the Euclidean ball centered at $x$ with radius $t$, and let
$|B(x,t)|$ denote its Lebesgue volume.  Define the Hardy--Littlewood maximal function
\begin{equation}\label{eq:maximal-function}
  Mh(x)
  :=\sup_{t>0}\frac1{|B(x,t)|}
  \int_{B(x,t)}|h(z)|\,dz.
\end{equation}
A pointwise Sobolev inequality of \citet{bojarski:hajlasz:1993} implies that,
for $g\in C^1(\mathbb R^d)$,
\begin{equation}\label{eq:sobolev-pointwise}
  |g(x)-g(y)|
  \leq C_d|x-y|
  \{M(|\nabla g|)(x)+M(|\nabla g|)(y)\}.
\end{equation}
Here $C_d<\infty$ depends only on the dimension.  Thus one may take
\begin{equation}\label{eq:sobolev-pointwise-G}
  G_g:=M(|\nabla g|),
\end{equation}
and define
\begin{equation}\label{eq:sobolev-class}
  \mathcal G_{\nabla,\pi}
  :=\left\{g\in C^1(\mathbb R^d)\cap L_0^2(\pi):
  M(|\nabla g|)\in L^4(\pi)\right\}.
\end{equation}

\begin{corollary}\label{thm:Sobolev-lift}
Under Assumption~\ref{ass:Wpsi2-envelope} and $ \Lambda\in L^4(\pi)$, every
$g\in\mathcal G_{\nabla,\pi}$ satisfies
\begin{equation}\label{eq:L2rate-sobolev}
  \|Q^kg\|_2
  \leq2C_d\|\Lambda\|_4
  \|M(|\nabla g|)\|_4r(k),
  \qquad k\geq1.
\end{equation}
Hence the Maxwell--Woodroofe condition holds, and the quenched CLT follows
under ergodicity.
\end{corollary}

The next proposition gives a sufficient condition for a $g$ belongs to $\mathcal G_{\nabla,\pi}$ when the invariant
density is bounded.

\begin{proposition}\label{prop:sobolev-sufficient}
Suppose $\pi(dx)=\varrho(x)\,dx$ and
$\varrho\in L^\infty(\mathbb R^d)$.  If
\[
  g\in C^1(\mathbb R^d)\cap L_0^2(\pi),
  \qquad
  \nabla g\in L^4(\real^d, dx),
\]
then $g\in\mathcal G_{\nabla,\pi}$.
\end{proposition}

\subsubsection{A polynomial-gradient construction}\label{poly_sufficient}

In this section, we  return directly to the line-segment bound and impose a pointwise growth condition on the gradient.
If, for some $K<\infty$ and $m\geq0$,
\begin{equation}\label{eq:polynomial-gradient-growth} |\nabla g(x)| \leq K(1+|x|^m) , \end{equation}
then there is a finite constant $C_{K,m}$ such that
\[
  |g(x)-g(y)|
  \leq C_{K,m}|x-y|
  \{(1+|x|^m)+(1+|y|^m)\}.
\]
Thus $G_g(x)=1+|x|^m$ is an increment envelope whenever it belongs to
$L^4(\pi)$.  This construction covers polynomial functions under invariant
laws with sufficiently many moments.

\begin{proposition}\label{prop:polynomial-gradient} 
Let $g\in C^1(\mathbb R^d)\cap L_0^2(\pi)$ satisfy \eqref{eq:polynomial-gradient-growth} for some $K<\infty$ and $m\geq0$. If \begin{equation}\label{eq:polynomial-moment}
\int_{\mathbb R^d}|x|^{4m}\,\pi(dx)<\infty,
\end{equation}
then $g$ admits an $L^4(\pi)$ increment envelope with $G_g(x)=1+|x|^m$. If, in addition, Assumption~\ref{ass:Wpsi2-envelope} holds and $\Lambda\in L^4(\pi)$, then 
\begin{equation}\label{eq:L2rate-polynomial} 
\|Q^kg\|_2 \leq 2C_{K,m}\|\Lambda\|_4 \|G_g\|_4\,r(k), \qquad k\geq1. \end{equation} Consequently, $g$ satisfies the Maxwell--Woodroofe condition. If the stationary chain is ergodic, then the quenched CLT \eqref{eq:quenched-clt} holds for $g$ and for $\pi$-almost every initial state. \end{proposition}

\section{Applications}\label{application}


\subsection{A nonlinear autoregressive chain}\label{nap}

 Let $a\in(0,1)$, let
$b:\mathbb R\to\mathbb R$ be measurable, and consider
\begin{equation}\label{nAR}
  X_{n+1}=F(X_n)+Z_{n+1},
  \qquad
  F(x):=ax+(1-a)b(x).
\end{equation}
 $(Z_n)_{n\geq1}$ are independent and identically distributed,
independent of $X_0$, with
$\mathbb E Z_1=0$ and $\mathbb E Z_1^2=\sigma_Z^2<\infty$.  Let $Q$ be the
transition kernel and use the Euclidean metric $\psi(x,y)=|x-y|$.

  Under a
 coupling, the same $Z_n$ is used for chains started at $x$
and $y$, so their one-step separation is $|F(x)-F(y)|$.  At the same time,   we use the quadratic Lyapunov function
\[
  \omega(x):=x^2.
\]
For $x\neq y$, define the local separation factor and the corresponding
Lyapunov ratio by
\begin{equation}\label{eq:NAR-zeta-kappa}
  \zeta(x,y):=\frac{|F(x)-F(y)|}{|x-y|},
  \qquad
  \kappa(x,y):=
  \frac{Q\omega(x)+Q\omega(y)+1}{\omega(x)+\omega(y)+1}.
\end{equation}
Because the $(Z_n)_{n\geq1}$  are centered,
\[
  Q\omega(x)=F(x)^2+\sigma_Z^2.
\]
Thus $\zeta(x,y)$ measures local expansion or contraction of the coupled
states, while $\kappa(x,y)$ measures the change in their combined Lyapunov
weight.

The following assumptions separate the requirements for the convergence bound
from the additional moments needed for the CLTs.
\begin{itemize}
  \item[$(N_1)$] There is a $\theta\in(0,1)$ such that
  \begin{equation}\label{eq:NAR-rho}
    \rho_\theta
    :=\sup_{x\neq y}\zeta(x,y)^\theta
      \kappa(x,y)^{1-\theta}<1.
  \end{equation}

  \item[$(N_2)$] Either $\sup_{x\neq y}\zeta(x,y)<\infty$, or $Q$ is weak
  Feller; that is, $Qf$ is continuous whenever $f$ is bounded and continuous.

  \item[$(N_3)$] The function $b$ is bounded and
  $\mathbb E|Z_1|^p<\infty$ for some $p>4$.
\end{itemize}

Condition $(N_1)$ is the essential drift--contraction condition: local
expansion in $F$ is allowed if it is compensated by sufficient decrease in the
Lyapunov factor.  Condition $(N_2)$ is a regularity condition used to obtain
existence and uniqueness of the invariant law.
Condition $(N_3)$ is not needed for the convergence bound; it supplies the higher
moments required for obtaining the CLT.  The
following consequence of \citet{qin:hobert:2022aihp} provides the
chain-side input for both applications in this section.

\begin{lemma}\label{GE_NAlinear}
Under $(N_1)$ and $(N_2)$, the chain has a unique invariant probability
measure $\pi$, and
\begin{equation}\label{NAlinear_rate}
  W_\psi(\delta_xQ^n,\pi)
  \leq \Lambda(x)\rho_\theta^n,
  \qquad n\geq1,
\end{equation}
where
\[
  \Lambda(x):=
  \frac{Q\omega(x)+\omega(x)+1}{1-\rho_\theta}.
\]
\end{lemma}

This convergence bound  does not require strict one-step Wasserstein contraction.  If
\begin{equation}\label{eq:NAR-no-GC}
  \sup_{x\neq y}\zeta(x,y)\geq1,
\end{equation}
then the following lemma shows that $\Delta_\psi(Q)\geq1$.

\begin{lemma}\label{not_GC}
Under~\eqref{eq:NAR-no-GC},
\[
  \sup_{x\neq y}
  \frac{W_\psi(\delta_xQ,\delta_yQ)}{|x-y|}\geq1.
\]
Thus the chain is not a strict one-step contraction in $W_\psi$.
\end{lemma}

For a concrete example, take $a=1/2$, $b(x)=-\sin x$, and $(Z_n)_{n\geq1}$ follow standard Gaussian.  Then $\sup_x|F'(x)|=1$, so strict one-step contraction fails.
Nevertheless, $(N_1)$ is verified in
\citet[Appendix~B]{qin:hobert:2022aihp}; $(N_2)$ follows from
$\zeta(x,y)\leq1$.  

The convergence rate alone is not yet enough to apply the Theorem~\ref{main_thm}: its
starting point envelope must be square integrable under $\pi$.  Since the
formula in Lemma~\ref{GE_NAlinear} grows quadratically in $x$, this reduces to
a fourth-moment calculation.   For
$q\geq1$, let
\[
  v_q(x):=|x|^q.
\]

\begin{lemma}\label{lem:NAR-moments}
Assume $(N_3)$.  For every $q\in[1,p]$, there is a finite constant $C_q$ such
that
\begin{equation}\label{eq:NAR-moment-bound}
  \sup_{n\geq0}Q^nv_q(x)
  \leq C_q(1+|x|^q),
  \qquad x\in\mathbb R.
\end{equation}
If $\delta_xQ^n\Rightarrow\pi$ for some invariant probability measure $\pi$,
then $\pi(v_q)<\infty$ for every $q\in[1,p]$.
\end{lemma}

Since $b$ is bounded, $\Lambda(x)\leq C(1+x^2)$.  Taking $q=4$ in
Lemma~\ref{lem:NAR-moments} gives $\Lambda\in L^2(\pi)$.  The geometric rate
$\rho_\theta^n$ automatically satisfies the Maxwell--Woodroofe summability
condition.  All assumptions of Theorem~\ref{main_thm} are now
verified, which yields the following CLT for Euclidean Lipschitz functions.

\begin{theorem}\label{NAR_clt}
Suppose $(N_1)$--$(N_3)$ hold.  Then every centered
Lipschitz function $g\in L^2(\pi)$ satisfies, for
$\pi$-almost every $x$,
\begin{equation}\label{eq:NAR-CLT}
  \frac{S_n(g)}{\sqrt n}
  \Rightarrow N(0,\sigma^2(g))
  \qquad\text{under }\mathbb P_x.
\end{equation}
\end{theorem}

\subsubsection{The quadratic function}\label{sec:AR-weighted}

Theorem~\ref{NAR_clt} does not cover $g(x) = x^2$, because this function
is not Euclidean Lipschitz.  This is precisely the type of situation for which
the weighted geometry was introduced.  We use
\[
  V(x):=1+x^2,
\]
which has the same  scale as the function of interest and is compatible with
the moment bounds in Lemma~\ref{lem:NAR-moments}.

 For $\alpha \in (0,1)$, define
\[
  \gamma_\alpha:=1-\frac{\alpha}{2},
  \qquad
  s_\alpha:=\frac{\gamma_\alpha}{1-\alpha}.
\]

We next show that  $g(x) = x^2$  belongs to the weighted Lipschitz class  induced by reweighting the Euclidean metric with $V$. 

\begin{lemma}\label{lem:x2-weighted-holder}
For every $A>0$ and $x,y\in\mathbb R$,
\begin{equation}\label{eq:x2-weighted-holder}
  |x^2-y^2|
  \leq C_\alpha|x-y|^\alpha
  \{A+V(x)+V(y)\}^{\gamma_\alpha},
  \qquad
  C_\alpha:=2^{1-\alpha/2}.
\end{equation}
Consequently, $g(x) = x^2/C_\alpha$ belongs to
$\mathcal G_{\psi,V}^{(\alpha,\gamma_\alpha,A)}$.
\end{lemma}

Proposition~\ref{prop:weighted-holder-lift} requires $V^{s_\alpha}\in L^2(\pi)$.  Since
$V(x)\asymp1+x^2$, this is equivalent to a moment of $|x|$ of order
$4s_\alpha$.  If $p>4$ is the moment condition in $(N_3)$, the available
moment bound is sufficient provided
\begin{equation}\label{eq:NAR-alpha-choice}
  0<\alpha<\frac{p-4}{p-2},
\end{equation}
which is equivalent to $4s_\alpha<p$. For such an $\alpha$, we can obtain a CLT for $g(x) = x^2$ by the following result. 

\begin{proposition}\label{prop:x2-CLT}
Suppose $(N_1)$--$(N_3)$ hold and choose $\alpha$ as in
\eqref{eq:NAR-alpha-choice}.  Let
$d_{\alpha,\gamma_\alpha,A}$ be the weighted path metric defined in
\eqref{eq:path-metric}.  There is a function
$\widetilde\Lambda_\alpha\in L^2(\pi)$ such that
\begin{equation}\label{eq:NAR-weighted-rate}
  W_{d_{\alpha,\gamma_\alpha,A}}(\delta_xQ^n,\pi)
  \leq\widetilde\Lambda_\alpha(x)\rho_\theta^{\alpha n}.
\end{equation}
Consequently, for
$g_2(x):=x^2-\pi(x^2)$ and for $\pi$-almost every $x$,
\begin{equation}\label{eq:NAR-x2-CLT}
  \frac{S_n(g_2)}{\sqrt n}
  \Rightarrow N(0,\sigma^2(g_2))
  \qquad\text{under }\mathbb P_x.
\end{equation}
\end{proposition}

\subsection{A subgeometric random walk on the $1$-dimensional torus}
\label{sec:circle_example}

We next give an example in which Theorem~\ref{thm:increment-envelope-lift} yields a CLT for a non-Lipschitz function of a random walk on the one-dimensional torus with a subgeometric convergence rate. We first use results from~\cite{Borda2021} to establish the required subgeometric convergence of the random walk. We then construct an increment envelope for a non-Lipschitz function and apply Theorem~\ref{thm:increment-envelope-lift} to establish its CLT, which is not covered by Theorem~1 of~\cite{Borda2021}.

Consider the $1$-dimensional torus $\mathbb R/\mathbb Z$ with the Euclidean metric and normalized Haar measure $\pi$. Given $x, y \in \real$, we define the Euclidean distance from the point $x + \mathbb Z$ to the point $y + \mathbb Z$ on the torus as $\rho(x,y) = \|x - y\|_{\mathbb R/\mathbb Z}:= \min_{m \in \mathbb Z} |x - y - m|$. 
Equivalently,
$
\rho(x,y)
:=
\min\{|x-y|,\,1-|x-y|\},
 x,y\in[0,1)$. Choose three points $\alpha_1,\alpha_2,\alpha_3\in \mathbb R/\mathbb Z$ forming a
badly approximable system in the sense that, for some $c>0$,
\begin{equation*}\label{eq:circle-badly-approx}
\max_{1\leq i\leq3}
\|m\alpha_i\|_{\mathbb R / \mathbb Z}
\geq c |m|^{-1/3},
\qquad
m\in\mathbb Z\setminus\{0\}. 
\end{equation*}
This is the condition~(6) of~\cite{Borda2021} specialized to $d=1$ and
$r=3$ for a badly approximable system.

Consider the random walk on $\mathbb R / \mathbb Z$
\begin{equation*}\label{eq:circle-chain}
X_{n+1}
=
X_n+Z_{n+1}\pmod 1,
\end{equation*}
where $Z_1,Z_2,\ldots$ are i.i.d. with distribution
$
\nu
:=
\frac1{6}
\sum_{i=1}^3
\left(
\delta_{\alpha_i}
+
\delta_{-\alpha_i}
\right)$.
Here $\delta_z$ denotes the Dirac probability measure concentrated at
$z $.  Theorems~6 and~7 of~\cite{Borda2021}, applied with $d=1$ and $r=3$,
give
\begin{equation}\label{eq:circle-W1-rate}
W_{\rho}(\delta_xQ^n,\pi)
\asymp
n^{-3/2},
\qquad
x\in  \mathbb R/\mathbb Z .
\end{equation}
Hence, $\pi$ is the unique invariant
probability measure and the stationary chain is ergodic.   Define
\begin{equation*}
\psi(x,y)
:=
\rho(x,y)^{1/2}.
\end{equation*}
Then, $\psi$ is a metric and the $2$-Wasserstein distance $W_{\psi,2}$ satisfies
\begin{equation*}\label{eq:circle-W2-identity}
W_{\psi,2}(\mu,\eta)^2
=
\inf_{\gamma\in\mathcal C(\mu,\eta)}
\int \rho(x,y)\,\gamma(dx,dy)
=
W_{\rho}(\mu,\eta).
\end{equation*}
Consequently, \eqref{eq:circle-W1-rate} gives
\begin{equation}\label{eq:circle-W2-rate}
W_{\psi,2}(\delta_xQ^n,\pi)
\asymp
n^{-3/4}.
\end{equation}
Thus, Assumption~\ref{ass:Wpsi2-envelope} holds with a  subgeometric rate
$r(n)=n^{-3/4}$. Consider the centered function
\begin{equation}\label{eq:circle-g}
g(x)
:=
\rho(x,0)^{1/3}
-
\pi\{\rho(\cdot,0)^{1/3}\}.
\end{equation}
Although $g\in L_0^2(\pi)$, it is not Lipschitz with respect to
$\psi=\rho^{1/2}$.  Indeed, for $t\downarrow0$,
\[
\frac{|g(t)-g(0)|}{\psi(t,0)}
=
\frac{t^{1/3}}{t^{1/2}}
=
t^{-1/6}
\longrightarrow\infty.
\]
Hence the $\psi$-Lipschitz condition required by Theorem~\ref{main_thm} does not apply.  Theorem~1 of~\cite{Borda2021} does not yield a CLT for this function
either.  

To state
the comparison with \citet{Borda2021} without conflicting with the rooted
$p$-Wasserstein distance in \eqref{eq:Wp}, define, for $0<q\leq1$,
\[
 \mathsf W_{\rho,q}(\mu,\eta)
 :=\inf_{\Gamma\in\mathcal C(\mu,\eta)}
 \int\rho(x,y)^q\,\Gamma(dx,dy).
\]
which is the Wasserstein distance defined in~\cite{Borda2021}.   The function $g$ is
$q$-H\"older with respect to $\rho$ exactly when $0<q\leq1/3$.  Theorems~6
and~7 of \citet{Borda2021} give
\[
 \mathsf W_{\rho,q}(\delta_xQ^n,\pi)^{1/q}\asymp n^{-3/2},
 \qquad
 \mathsf W_{\rho,q}(\delta_xQ^n,\pi)\asymp n^{-3q/2}.
\]
The summability assumption in Theorem~1 of that paper therefore requires
$q>2/3$, which is incompatible with $q\leq1/3$.  
Consequently, the conditions of Theorem~1 of~\cite{Borda2021}
are not satisfied.

By contrast, all conditions of Theorem~\ref{thm:increment-envelope-lift} can be verified. We next identify an increment envelope for this function $g$. 

\begin{proposition}\label{prop:circle-increment-envelope}
Let $g$ be defined by~\eqref{eq:circle-g}, and define
\begin{equation*}\label{eq:circle-G}
G_g(x)
:=
\rho(x,0)^{-1/6}
\mathbf 1_{\{x\neq0\}}.
\end{equation*}
Then $G_g\in L^4(\pi)$ and
\begin{equation*}\label{eq:circle-increment}
|g(x)-g(y)|
\leq
\psi(x,y)\{G_g(x)+G_g(y)\},
\qquad x,y\in\mathbb R/\mathbb Z .
\end{equation*}
\end{proposition}

Since Assumption~\ref{ass:Wpsi2-envelope} holds with
$\Lambda\equiv C$ and $r(n)=n^{-3/4}$, we have
$\Lambda\in L^4(\pi)$.  Proposition~\ref{prop:circle-increment-envelope}
also gives $G_g\in L^4(\pi)$.  Taking $a=b=2$ in
Theorem~\ref{thm:increment-envelope-lift} therefore yields
\begin{equation*}\label{eq:circle-L2-rate}
\|Q^ng\|
\lesssim
n^{-3/4}.
\end{equation*}
Since
\[
\sum_{n=1}^{\infty}
\frac{\|Q^ng\|}{\sqrt n}
\lesssim
\sum_{n=1}^{\infty}n^{-5/4}
<\infty,
\]
the required projective summability condition holds, and hence so does the Maxwell--Woodroofe criterion.  Consequently, for
$\pi$--almost every $x$,
\begin{equation*}\label{eq:circle-clt}
\frac1{\sqrt n}
\sum_{k=0}^{n-1}g(X_k)
\Longrightarrow
N(0,\sigma_g^2)
\end{equation*}
under $\mathbb P_x$.

Finally, this example lies genuinely outside a total variation
convergence framework.  For every fixed $n$, the measure
$\delta_xQ^n$ has finite support, whereas $\pi$ is nonatomic.
Hence $\delta_xQ^n$ and $\pi$ are mutually singular for every
finite $n$.  In particular, their total variation distance is
maximal:
\[
\|\delta_xQ^n-\pi\|_{\mathrm{TV}}=1,
\qquad
x\in\mathbb R/\mathbb Z,\quad n\geq1.
\]

\appendix

\section{Proofs for Section~\ref{main_result}}\label{app:section3}

The proofs follow the structure emphasized in the main text.  

\subsection{Proof of Theorem~\ref{main_thm}}
\begin{proof}
The proof has two steps.  We first sum the $L^2$ bias bounds to control
$V_ng$, and then reverse the order of summation in the Maxwell--Woodroofe
series.  Let $g\in\mathcal G_\psi\cap L_0^2(\pi)$.  By
\eqref{eq:key-L2-decay-norm},
\[
  \|Q^kg\|_2\leq\|\Lambda\|_2r(k),
  \qquad k\geq1.
\]
The term $Q^0g=g$ is kept separate.  Hence
\[
  \left\|\sum_{k=0}^{n-1}Q^kg\right\|_2
  \leq\|g\|_2+\|\Lambda\|_2\sum_{k=1}^{n-1}r(k).
\]
Tonelli's theorem gives
\[
\begin{aligned}
  \sum_{n=1}^{\infty}n^{-3/2}
  \left\|\sum_{k=0}^{n-1}Q^kg\right\|_2
  &\leq\|g\|_2\sum_{n=1}^{\infty}n^{-3/2}\\
  &\quad+\|\Lambda\|_2
  \sum_{k=1}^{\infty}r(k)
  \sum_{n=k+1}^{\infty}n^{-3/2}.
\end{aligned}
\]
Since
\[
  \sum_{n=k+1}^{\infty}n^{-3/2}
  \leq\int_k^\infty t^{-3/2}\,dt
  =\frac{2}{\sqrt k},
\]
condition~\eqref{rate_fun} proves~\eqref{suff_cond_clt}.  The quenched CLT
follows from Theorem~\ref{clt_anypoint} under ergodicity.  If
$0<\operatorname{Lip}_\psi(g)<\infty$, apply the result to
$g/\operatorname{Lip}_\psi(g)$ and rescale.  
\end{proof}

\subsection{Proof of Corollary~\ref{check:solvePE}}
\begin{proof}
Here the stronger summability of $r$ makes the series of $L^2$ norms directly
summable.  Equation~\eqref{eq:key-L2-decay-norm} implies
\[
  \sum_{n=0}^{\infty}\|Q^ng\|_2
  \leq\|g\|_2+\|\Lambda\|_2\sum_{n=1}^{\infty}r(n)<\infty.
\]
Let $h_N=\sum_{n=0}^NQ^ng$.  Then $h_N\to h$ in $L^2(\pi)$ for some $h$.
Since $Q$ is an $L^2(\pi)$ contraction, $Qh_N\to Qh$.  Moreover,
\[
  h_N-Qh_N=g-Q^{N+1}g\longrightarrow g
  \qquad\text{in }L^2(\pi),
\]
so $h-Qh=g$.
\end{proof}

\subsection{Proof of Proposition~\ref{MA_GC}}
\begin{proof}
  Write
$\Delta=\Delta_\psi(Q)<1$.   Let $h\in\mathcal G_\psi$.  Subtracting a constant if necessary, assume
$h(x_0)=0$.  Then
\[
  |Qh(x)-Qh(y)|
  \leq W_\psi(\delta_xQ,\delta_yQ)
  \leq\Delta\psi(x,y),
\]
so $Qh$ is $\Delta$-Lipschitz.  Kantorovich--Rubinstein duality yields
\[
  W_\psi(\mu Q,\nu Q)
  \leq\Delta W_\psi(\mu,\nu),
  \qquad \mu,\nu\in\mathcal P_\psi^1.
\]
Using $\pi Q=\pi$ and iterating,
\[
  W_\psi(\delta_xQ^n,\pi)
  \leq\Delta^nW_\psi(\delta_x,\pi).
\]
Thus~\eqref{w_subge} holds with
$r(n)=\Delta^n$ and
\[
  \Lambda(x)=W_\psi(\delta_x,\pi)
  =\int\psi(x,y)\,\pi(dy).
\]
Choose $x_0$ so that $\pi\{\psi(x_0,\cdot)^2\}<\infty$.  The triangle
inequality gives
\[
  \Lambda(x)
  \leq\psi(x,x_0)+\pi\{\psi(x_0,\cdot)\},
\]
so $\Lambda\in L^2(\pi)$.  The geometric sequence satisfies both
\eqref{rate_fun} and~\eqref{eq:poisson-rate}.  The conclusions follow from
Theorem~\ref{main_thm} and Corollary~\ref{check:solvePE}.
\end{proof}

\subsection{Proof of Proposition~\ref{lem:path-metric-geometry}}
\begin{proof}
 The upper
bound in~\eqref{eq:path-lower-upper} follows from the one-edge path
$(x,y)$.  For the lower bound, consider an arbitrary finite path
$x=z_0,\ldots,z_m=y$.  Since
$A+V(z_{i-1})+V(z_i)\geq A$,
\[
\begin{aligned}
  \sum_{i=1}^m c_{\alpha,A}(z_{i-1},z_i)
  &\geq A^{1-\alpha}\sum_{i=1}^m
  \psi(z_{i-1},z_i)^\alpha\\
  &\geq A^{1-\alpha}
  \left\{\sum_{i=1}^m\psi(z_{i-1},z_i)\right\}^\alpha\\
  &\geq A^{1-\alpha}\psi(x,y)^\alpha.
\end{aligned}
\]
The second inequality uses subadditivity of $t\mapsto t^\alpha$, and the
third uses the triangle inequality for $\psi$.  Taking the infimum over paths
proves~\eqref{eq:path-lower-upper}.

Symmetry is immediate, the lower bound gives nondegeneracy, and concatenation
of paths gives the triangle inequality.  Hence $d_{\alpha,A}$ is a metric.  If
$d_{\alpha,A}(x_n,x)\to0$, the lower bound implies $\psi(x_n,x)\to0$.
Conversely, if $\psi(x_n,x)\to0$, continuity of $V$ and the upper bound
$d_{\alpha,A}\leq c_{\alpha,A}$ imply $d_{\alpha,A}(x_n,x)\to0$.  Thus the
two metrics generate the same topology and Borel sets.

To verify completeness, let $(x_n)$ be $d_{\alpha,A}$-Cauchy.  The lower
bound shows that it is $\psi$-Cauchy, so $x_n\to x$ in $\psi$ for some $x$.
The upper bound and continuity of $V$ then give $d_{\alpha,A}(x_n,x)\to0$.
Separability follows because the two metrics generate the same topology.

It remains to identify the Lipschitz class.  If
$\operatorname{Lip}_{d_{\alpha,A}}(g)\leq1$, then
\[
  |g(x)-g(y)|
  \leq d_{\alpha,A}(x,y)
  \leq c_{\alpha,A}(x,y),
\]
so $g\in\mathcal G_{\psi,V}^{(\alpha,A)}$.  Conversely, if $g$ is in that
class, then along any path $x=z_0,\ldots,z_m=y$,
\[
  |g(x)-g(y)|
  \leq\sum_{i=1}^m c_{\alpha,A}(z_{i-1},z_i).
\]
Taking the infimum over all paths yields
$|g(x)-g(y)|\leq d_{\alpha,A}(x,y)$.
\end{proof}

\subsection{Proof of Theorem~\ref{thm:direct-lift}}
\begin{proof}
 Fix a coupling
$\Gamma\in\mathcal C(\mu,\nu)$.  Since
$d_{\alpha,A}\leq c_{\alpha,A}$, H\"older's inequality gives
\[
\begin{aligned}
  \int d_{\alpha,A}\,d\Gamma
  &\leq\int\psi(x,y)^\alpha
  \{A+V(x)+V(y)\}^{1-\alpha}\,d\Gamma\\
  &\leq
  \left\{\int\psi(x,y)\,d\Gamma\right\}^\alpha
  \left\{A+\mu(V)+\nu(V)\right\}^{1-\alpha}.
\end{aligned}
\]
Taking the infimum over $\Gamma$ proves~\eqref{eq:direct-lift}.
\end{proof}

\subsection{Proof of Proposition~\ref{prop:direct-lift}}
\begin{proof}
 Apply Theorem~\ref{thm:direct-lift} with
$\mu=\delta_xQ^n$ and $\nu=\pi$.  Since $V\in L^2(\pi)$,
$\pi(V)<\infty$.  Assumption~\eqref{eq:direct-moment-iterate} gives
\[
  A+Q^nV(x)+\pi(V)\leq C_0\{1+V(x)\}
\]
for a finite constant $C_0$.  Together with~\eqref{w_subge}, this yields
\[
  W_{d_{\alpha,A}}(\delta_xQ^n,\pi)
  \leq C\Lambda(x)^\alpha
  \{1+V(x)\}^{1-\alpha}r(n)^\alpha,
\]
which is~\eqref{eq:direct-transferred-envelope}.  H\"older's inequality gives
\[
  \pi(\widetilde\Lambda_\alpha^2)
  \leq C^2\{\pi(\Lambda^2)\}^\alpha
  \{\pi[(1+V)^2]\}^{1-\alpha}<\infty.
\]
By Proposition~\ref{lem:path-metric-geometry},
$(\mathcal X,d_{\alpha,A})$ is Polish and its unit Lipschitz class is
$\mathcal G_{\psi,V}^{(\alpha,A)}$.  Theorem~\ref{main_thm} therefore applies
with rate $r(n)^\alpha$.  This proves the Maxwell--Woodroofe condition and,
under ergodicity, the quenched CLT.
\end{proof}

\subsection{Proof of Theorem~\ref{thm:weighted-holder-lift}}
\begin{proof}
The geometric part is unchanged from the canonical weighted construction.  In
the measure comparison, the identity $\gamma=s(1-\alpha)$ is what allows
H\"older's inequality to isolate the $s$th moment of $V$.  The geometric
assertions follow by repeating the proof of
Proposition~\ref{lem:path-metric-geometry} with
$c_{\alpha,A}$ and $d_{\alpha,A}$ replaced by
$c_{\alpha,\gamma,A}$ and $d_{\alpha,\gamma,A}$.  The lower bound is now
$A^\gamma\psi(x,y)^\alpha$, and the remaining topology, completeness, and
Lipschitz-class arguments are unchanged.

It remains to prove~\eqref{eq:wh-transfer}.  Let
$s=\gamma/(1-\alpha)$ and fix a coupling
$\Gamma\in\mathcal C(\mu,\nu)$.  Since
$d_{\alpha,\gamma,A}\leq c_{\alpha,\gamma,A}$ and
$\gamma=s(1-\alpha)$, H\"older's inequality gives
\[
\begin{aligned}
  \int d_{\alpha,\gamma,A}\,d\Gamma
  &\leq\int\psi(x,y)^\alpha
  \{A+V(x)+V(y)\}^{s(1-\alpha)}\,d\Gamma\\
  &\leq\left\{\int\psi(x,y)\,d\Gamma\right\}^\alpha
  \left\{\int\{A+V(x)+V(y)\}^s\,d\Gamma\right\}^{1-\alpha}.
\end{aligned}
\]
For $u,v,w\geq0$,
\[
  (u+v+w)^s\leq c_s(u^s+v^s+w^s),
\]
with $c_s$ as in the theorem.  Therefore
\[
  \int\{A+V(x)+V(y)\}^s\,d\Gamma
  \leq c_s\{A^s+\mu(V^s)+\nu(V^s)\}.
\]
Taking the infimum over $\Gamma$ proves~\eqref{eq:wh-transfer}.
\end{proof}

\subsection{Proof of Proposition~\ref{prop:weighted-holder-lift}}
\begin{proof}
 Apply
Theorem~\ref{thm:weighted-holder-lift} with
$\mu=\delta_xQ^n$ and $\nu=\pi$.  Since $V^s\in L^2(\pi)$,
$\pi(V^s)<\infty$.  Assumption~\eqref{eq:weighted-moment-iterate} gives
\[
  A^s+Q^nV^s(x)+\pi(V^s)
  \leq C_0\{1+V(x)^s\}
\]
for a finite constant $C_0$.  Together with~\eqref{w_subge}, this yields
\[
  W_{d_{\alpha,\gamma,A}}(\delta_xQ^n,\pi)
  \leq C\Lambda(x)^\alpha
  \{1+V(x)^s\}^{1-\alpha}r(n)^\alpha.
\]
Thus~\eqref{eq:direct-transferred-envelope-general} holds with the stated
$\widetilde\Lambda_{\alpha,\gamma}$.  H\"older's inequality gives
\[
  \pi(\widetilde\Lambda_{\alpha,\gamma}^2)
  \leq C^2\{\pi(\Lambda^2)\}^\alpha
  \{\pi[(1+V^s)^2]\}^{1-\alpha}<\infty.
\]
By Theorem~\ref{thm:weighted-holder-lift},
$(\mathcal X,d_{\alpha,\gamma,A})$ is Polish and its unit Lipschitz class is
$\mathcal G_{\psi,V}^{(\alpha,\gamma,A)}$.  Theorem~\ref{main_thm} therefore
applies with rate $r(n)^\alpha$.  This proves the Maxwell--Woodroofe condition
and, under ergodicity, the quenched CLT.
\end{proof}

\subsection{Proof of Theorem~\ref{thm:increment-envelope-lift}}
\begin{proof}
  Fix
$k\geq1$.  Since $G_g\in L^{2b}(\pi)\subset L^2(\pi)$ and $\pi$ is
invariant,
\[
  \int Q^kG_g^2(x)\,\pi(dx)=\pi(G_g^2)<\infty.
\]
Hence $Q^kG_g^2(x)<\infty$ for $\pi$-almost every $x$.  Fix such an $x$ and
let $\Gamma_{k,x}$ be an optimal $W_{\psi, 2}$ coupling of
$\delta_xQ^k$ and $\pi$.  The increment inequality and Cauchy--Schwarz give
\[
\begin{aligned}
  |Q^kg(x)-\pi(g)|^2
  &\leq C_g^2
  \left(\int \psi(y, z)^2\,d\Gamma_{k,x}\right)
  \left(\int\{G_g(y)+G_g(z)\}^2\,d\Gamma_{k,x}\right)\\
  &\leq2C_g^2W_{\psi, 2}(\delta_xQ^k,\pi)^2
  \{Q^kG_g^2(x)+\pi(G_g^2)\}.
\end{aligned}
\]
  By
Assumption~\ref{ass:Wpsi2-envelope},
\[
  |Q^kg(x)|^2
  \leq2C_g^2\Lambda(x)^2r(k)^2
  \{Q^kG_g^2(x)+\pi(G_g^2)\}.
\]
Integrating with respect to $\pi$ yields
\[
  \|Q^kg\|_2^2
  \leq2C_g^2r(k)^2(I_1+I_2),
\]
where
\[
  I_1:=\pi(\Lambda^2Q^kG_g^2),
  \qquad
  I_2:=\pi(\Lambda^2)\pi(G_g^2).
\]
We next bound $I_1$ and $I_2$ separately.  For $I_1$, H\"older's
inequality with conjugate exponents $a$ and $b$ gives
\begin{align}
  I_1
  &=\bpi(\Lambda^2\bQ^kG_g^2)\notag\\
  &\leq
  \bpi(\Lambda^{2a})^{1/a}
  \bpi\{(\bQ^kG_g^2)^b\}^{1/b}.
  \label{eq:I1-holder}
\end{align}
For each $x$, $\bQ^k(x,\cdot)$ is a probability measure and
$t\mapsto t^b$ is convex because $b>1$.  Jensen's inequality therefore yields
\begin{equation}\label{eq:I1-jensen}
  \{\bQ^kG_g^2(x)\}^b
  =\left\{\int G_g(y)^2\,\bQ^k(x,dy)\right\}^b
  \leq\int G_g(y)^{2b}\,\bQ^k(x,dy)
  =\bQ^kG_g^{2b}(x).
\end{equation}
Substituting \eqref{eq:I1-jensen} into \eqref{eq:I1-holder}, and then using
invariance of $\bpi$, gives
\begin{align}
  I_1
  &\leq
  \bpi(\Lambda^{2a})^{1/a}
  \bpi(\bQ^kG_g^{2b})^{1/b}\notag\\
  &=\bpi(\Lambda^{2a})^{1/a}
    \bpi(G_g^{2b})^{1/b}\notag\\
  &=\|\Lambda\|_{2a}^2\|G_g\|_{2b}^2.
  \label{eq:I1-final}
\end{align}
Here the equality in the second line uses
$\bpi\bQ^k=\bpi$.

For $I_2$, apply H\"older's inequality separately to each factor, using the
constant function $1$.  Since $\bpi$ is a probability measure,
\begin{align*}
  \bpi(\Lambda^2)
  &=\bpi(\Lambda^2\cdot1)
    \leq\bpi(\Lambda^{2a})^{1/a}\bpi(1^b)^{1/b}
    =\|\Lambda\|_{2a}^2,\\
  \bpi(G_g^2)
  &=\bpi(G_g^2\cdot1)
    \leq\bpi(G_g^{2b})^{1/b}\bpi(1^a)^{1/a}
    =\|G_g\|_{2b}^2.
\end{align*}
Consequently,
\begin{equation}\label{eq:I2-final}
  I_2
  =\bpi(\Lambda^2)\bpi(G_g^2)
  \leq\|\Lambda\|_{2a}^2\|G_g\|_{2b}^2.
\end{equation}
Combining \eqref{eq:I1-final} and \eqref{eq:I2-final},
\[
  I_1+I_2
  \leq2\|\Lambda\|_{2a}^2\|G_g\|_{2b}^2.
\]
Therefore,
\[
  \|\bQ^kg\|_2^2
  \leq4C_g^2\|\Lambda\|_{2a}^2\|G_g\|_{2b}^2r(k)^2,
\]
and taking square roots gives
\[
  \|\bQ^kg\|_2
  \leq2C_g\|\Lambda\|_{2a}\|G_g\|_{2b}r(k).
\]
The summation argument in Theorem~\ref{main_thm}, with this bound in place of
\eqref{eq:key-L2-decay-norm}, proves the Maxwell--Woodroofe condition.
Theorem~\ref{clt_anypoint} gives the quenched CLT under ergodicity.
\end{proof}

\subsection{Proof of Corollary~\ref{thm:Sobolev-lift}}
\begin{proof}
Equation~\eqref{eq:sobolev-pointwise} gives the increment envelope with
$C_g=C_d$ and $G_g=M(|\nabla g|)$.  Apply
Theorem~\ref{thm:increment-envelope-lift}.
\end{proof}

\subsection{Proof of Proposition~\ref{prop:sobolev-sufficient}}
\begin{proof}
The strong $(4,4)$ Hardy--Littlewood maximal inequality gives
\[
\begin{aligned}
  \|M(|\nabla g|)\|_4^4
  &=\int M(|\nabla g|)(x)^4\varrho(x)\,dx\\
  &\leq\|\varrho\|_\infty
  \|M(|\nabla g|)\|_{L^4(dx)}^4\\
  &\leq\|\varrho\|_\infty C_{d,4}^4
  \|\nabla g\|_{L^4(dx)}^4.
\end{aligned}
\]
Taking fourth roots proves the result.
\end{proof}

\subsection{Proof of Proposition~\ref{prop:polynomial-gradient}}
\begin{proof}
For $x,y\in\mathbb R^d$, the fundamental theorem of calculus and
\eqref{eq:polynomial-gradient-growth} give
\[
 |g(x)-g(y)|
 \leq K|x-y|\int_0^1
 \{1+|(1-t)y+tx|^m\}\,dt.
\]
For $m\geq0$,
\[
 |(1-t)y+tx|^m
 \leq c_m\{|x|^m+|y|^m\},
 \qquad c_m:=\max(1,2^{m-1}).
\]
Consequently,
\[
 |g(x)-g(y)|
 \leq C_{K,m}|x-y|
 \{G_g(x)+G_g(y)\},
 \qquad G_g(x):=1+|x|^m,
\]
where one may take $C_{K,m}=K\max(1,c_m)$.  Condition
\eqref{eq:polynomial-moment} implies $G_g\in L^4(\pi)$.  The final
$L^2(\pi)$ bound and the CLT now follow from
Theorem~\ref{thm:increment-envelope-lift} with $a=b=2$.
\end{proof}

\section{Proofs for Section~\ref{application}}\label{app:applications}

\subsection{Proof of Lemma~\ref{GE_NAlinear}}
\begin{proof}
We apply \citet[Theorem~2.5]{qin:hobert:2022aihp} with
$\omega(x)=x^2$ and $\psi(x,y)=|x-y|$.  Its growth condition holds because
\[
  |x-y|\leq|x|+|y|
  \leq x^2+y^2+\tfrac12
  \leq\omega(x)+\omega(y)+1.
\]
Under this coupling, both chains use the same $Z_n$, so
\[
  W_\psi(\delta_xQ,\delta_yQ)
  \leq|F(x)-F(y)|
  =\zeta(x,y)\psi(x,y).
\]
Condition $(N_1)$ is the joint drift--contraction condition in that theorem,
with contraction factor $\rho_\theta$.  It follows that an invariant
probability measure exists and that~\eqref{NAlinear_rate} holds.
Condition $(N_2)$ gives uniqueness through
\citet[Proposition~2.7]{qin:hobert:2022aihp}.
\end{proof}

\subsection{Proof of Lemma~\ref{lem:NAR-moments}}
\begin{proof}
Let $B:=\sup_x|b(x)|<\infty$.
Iterating~\eqref{nAR} gives
\[
  X_n
  =a^nx
  +(1-a)\sum_{j=0}^{n-1}a^{n-1-j}b(X_j)
  +\sum_{j=1}^na^{n-j}Z_j.
\]
Fix $q\in[1,p]$.  Minkowski's inequality yields
\[
\begin{aligned}
  \{\mathbb E_x|X_n|^q\}^{1/q}
  &\leq a^n|x|
  +(1-a)B\sum_{j=0}^{n-1}a^{n-1-j}
  +\|Z_1\|_{L^q}\sum_{j=1}^na^{n-j}\\
  &\leq |x|+B+\frac{\|Z_1\|_{L^q}}{1-a}.
\end{aligned}
\]
Raising both sides to the $q$th power and absorbing constants proves
\eqref{eq:NAR-moment-bound}.

If $\delta_xQ^n\Rightarrow\pi$, then $v_q(z)=|z|^q$ is nonnegative and lower
semicontinuous.  The Portmanteau theorem gives
\[
  \pi(v_q)\leq\liminf_{n\to\infty}Q^nv_q(x)<\infty.
\]
\end{proof}

\subsection{Proof of Theorem~\ref{NAR_clt}}
\begin{proof}
  Lemma~\ref{GE_NAlinear} gives the invariant law and the geometric Wasserstein
bound~\eqref{NAlinear_rate}.  Uniqueness of the invariant probability measure
implies ergodicity of the stationary chain.  Since $b$ is bounded,
\[
  Q\omega(x)+\omega(x)+1
  =F(x)^2+\sigma_Z^2+x^2+1
  \leq C(1+x^2).
\]
Thus $\Lambda(x)\leq C(1+x^2)$.  Lemma~\ref{lem:NAR-moments}, with $q=4$,
gives $\pi(x^4)<\infty$, and hence $\Lambda\in L^2(\pi)$.  The rate
$r(n)=\rho_\theta^n$ satisfies~\eqref{rate_fun}.  The result follows from
Theorem~\ref{main_thm}.
\end{proof}

\subsection{Proof of Lemma~\ref{not_GC}}
\begin{proof}
The identity function $\iota(z):=z$ is $1$-Lipschitz under the Euclidean
metric.  Kantorovich--Rubinstein duality gives
\[
  W_\psi(\delta_xQ,\delta_yQ)
  \geq|Q\iota(x)-Q\iota(y)|
  =|F(x)-F(y)|.
\]
Divide by $|x-y|$ and take the supremum.
\end{proof}

\subsection{Proof of Lemma~\ref{lem:x2-weighted-holder}}
\begin{proof}
 Let $\delta:=|x-y|$ and $R:=(x^2+y^2)^{1/2}$.  Since
$|x+y|\leq\sqrt2R$ and $\delta\leq\sqrt2R$,
\[
\begin{aligned}
  |x^2-y^2|
  &=\delta|x+y|\\
  &\leq\sqrt2\,\delta R\\
  &=\sqrt2\,\delta^\alpha\delta^{1-\alpha}R\\
  &\leq2^{1-\alpha/2}\delta^\alpha R^{2-\alpha}.
\end{aligned}
\]
Since $A+V(x)+V(y)\geq x^2+y^2=R^2$ and
$\gamma_\alpha=1-\alpha/2$, the last expression is bounded by the
right-hand side of~\eqref{eq:x2-weighted-holder}.
\end{proof}

\subsection{Proof of Proposition~\ref{prop:x2-CLT}}
\begin{proof}
 By~\eqref{eq:NAR-alpha-choice},
\[
  4s_\alpha=\frac{2(2-\alpha)}{1-\alpha}<p.
\]
Lemma~\ref{lem:NAR-moments} therefore gives
$\pi(V^{2s_\alpha})<\infty$, so $V^{s_\alpha}\in L^2(\pi)$.  The same lemma,
used at moment order $2s_\alpha$, yields
\[
  \sup_{n\geq1}Q^nV^{s_\alpha}(x)
  \leq C\{1+V(x)^{s_\alpha}\}.
\]
The envelope in~\eqref{NAlinear_rate} belongs to $L^2(\pi)$ by the proof of
Theorem~\ref{NAR_clt}.  Proposition~\ref{prop:weighted-holder-lift}, with
$r(n)=\rho_\theta^n$, gives~\eqref{eq:NAR-weighted-rate} and the quenched CLT
for every centered function in
$\mathcal G_{\psi,V}^{(\alpha,\gamma_\alpha,A)}$.  Lemma
\ref{lem:x2-weighted-holder} and rescaling give the stated result for $g_2$.
\end{proof}

\subsection{Proof of Proposition~\ref{prop:circle-increment-envelope}}
\label{app:circle-envelope}

\begin{proof}
We first establish the elementary inequality
\begin{equation}\label{eq:circle-scalar-ineq}
|u^{1/3}-v^{1/3}|
\leq
|u-v|^{1/2}
\left\{
u^{-1/6}\mathbf 1_{\{u>0\}}
+
v^{-1/6}\mathbf 1_{\{v>0\}}
\right\},
\qquad u,v\geq0.
\end{equation}
By symmetry, suppose $u\geq v$.  If $v=0$, then
\[
u^{1/3}
=
u^{1/2}u^{-1/6},
\]
so~\eqref{eq:circle-scalar-ineq} holds.  If $v>0$, then
\[
u^{1/3}-v^{1/3}
=
\frac{u-v}
{u^{2/3}+u^{1/3}v^{1/3}+v^{2/3}}.
\]
Since $u\geq v$, it follows that
\[
u^{1/3}-v^{1/3}
\leq
(u-v)^{1/2}v^{-1/6},
\]
and~\eqref{eq:circle-scalar-ineq} follows.

Now take
\[
u=\rho(x,0),
\qquad
v=\rho(y,0).
\]
Since the distance from a fixed point is $1$--Lipschitz,
\[
|\rho(x,0)-\rho(y,0)|
\leq
\rho(x,y).
\]
The centering constant in~\eqref{eq:circle-g} cancels from
$g(x)-g(y)$, so~\eqref{eq:circle-scalar-ineq} gives
\[
|g(x)-g(y)|
\leq
\rho(x,y)^{1/2}
\{G_g(x)+G_g(y)\}
=
\psi(x,y)\{G_g(x)+G_g(y)\}.
\]

It remains to verify the required integrability.  Identifying
$\mathbb R/\mathbb Z$ with $[0,1)$ and using normalized Haar measure,
\[
\|G_g\|_{L^4(\pi)}^4
=
2\int_0^{1/2}t^{-2/3}\,dt
<\infty.
\]
Thus $G_g\in L^4(\pi)$.
\end{proof}


\begin{thebibliography}{21}

\bibitem[{Bojarski and Haj{\l}asz(1993)}]{bojarski:hajlasz:1993}
Bojarski, B. and Haj{\l}asz, P. (1993).
\newblock Pointwise inequalities for Sobolev functions and some applications.
\newblock \textit{Studia Mathematica, 106}(1), 77--92.

\bibitem[{Borda(2021)}]{Borda2021}
Borda, B. (2021).
\newblock Equidistribution of random walks on compact groups II.\ The Wasserstein metric.
\newblock \textit{Bernoulli, 27}(4), 2598--2623.

\bibitem[{Butkovsky(2014)}]{butkovsky:2014}
Butkovsky, O. (2014).
\newblock Subgeometric rates of convergence of Markov processes in the Wasserstein metric.
\newblock \textit{The Annals of Applied Probability, 24}(2), 526--552.

\bibitem[{Chan and Geyer(1994)}]{chan:geyer:1994}
Chan, K.~S. and Geyer, C.~J. (1994).
\newblock Discussion: Markov chains for exploring posterior distributions.
\newblock \textit{The Annals of Statistics, 22}(4), 1747--1758.

\bibitem[{Cuny and Merlev{\`e}de(2014)}]{cuny:merlevede:2014}
Cuny, C. and Merlev{\`e}de, F. (2014).
\newblock On martingale approximations and the quenched weak invariance principle.
\newblock \textit{The Annals of Probability, 42}(2), 760--793.

\bibitem[{Czapla et~al.(2024)}]{czapla:horbacz:wojewodka:2024}
Czapla, D., Horbacz, K. and Wojewodka-Sciazko, H. (2024).
\newblock The central limit theorem for Markov processes that are exponentially ergodic in the bounded-Lipschitz norm.
\newblock \textit{Qualitative Theory of Dynamical Systems, 23}, 7.

\bibitem[{Durmus and Moulines(2015)}]{durmus:moulines:2015}
Durmus, A. and Moulines, E. (2015).
\newblock Quantitative bounds of convergence for geometrically ergodic Markov chain in the Wasserstein distance with application to the Metropolis adjusted Langevin algorithm.
\newblock \textit{Statistics and Computing, 25}(1), 5--19.

\bibitem[{Eberle and Majka(2019)}]{eberle:majka:2019}
Eberle, A. and Majka, M.~B. (2019).
\newblock Quantitative contraction rates for Markov chains on general state spaces.
\newblock \textit{Electronic Journal of Probability, 24}, Paper~26, 1--36.

\bibitem[{Gulgowski et~al.(2019)Gulgowski, Hille, Szarek and Ziemla{\'n}ska}]{gulgowski:hille:szarek:ziemlanska:2019}
Gulgowski, J., Hille, S.~C., Szarek, T. and Ziemla{\'n}ska, M.~A. (2019).
\newblock Central limit theorem for some non-stationary Markov chains.
\newblock \textit{Studia Mathematica, 246}(2), 109--131.

\bibitem[{Hairer, Mattingly and Scheutzow(2011)}]{hairer:mattingly:scheutzow:2011}
Hairer, M., Mattingly, J.~C. and Scheutzow, M. (2011).
\newblock Asymptotic coupling and a general form of Harris' theorem with applications to stochastic delay equations.
\newblock \textit{Probability Theory and Related Fields, 149}, 223--259.

\bibitem[{Jarner and Roberts(2002)}]{jarner:robert:2002}
Jarner, S.~F. and Roberts, G.~O. (2002).
\newblock Polynomial convergence rates of Markov chains.
\newblock \textit{The Annals of Applied Probability, 12}(1), 224--247.

\bibitem[{Jones(2004)}]{jones:2004}
Jones, G.~L. (2004).
\newblock On the Markov chain central limit theorem.
\newblock \textit{Probability Surveys, 1}, 299--320.

\bibitem[{Komorowski and Walczuk(2012)}]{komorowski:walczuk:2012}
Komorowski, T. and Walczuk, A. (2012).
\newblock Central limit theorem for Markov processes with spectral gap in the Wasserstein metric.
\newblock \textit{Stochastic Processes and their Applications, 122}(5), 2155--2184.

\bibitem[{Kulik(2018)}]{kulik:2018}
Kulik, A. (2018).
\newblock \textit{Ergodic Behavior of Markov Processes: With Applications to Limit Theorems}.
\newblock De Gruyter Studies in Mathematics, Vol.~67. De Gruyter, Berlin/Boston.

\bibitem[{Maxwell and Woodroofe(2000)}]{maxwell:woodroofe:2000}
Maxwell, M. and Woodroofe, M. (2000).
\newblock Central limit theorems for additive functionals of Markov chains.
\newblock \textit{The Annals of Probability, 28}(2), 713--724.

\bibitem[{Qin and Hobert(2019)}]{qin:hobe:2019}
Qin, Q. and Hobert, J.~P. (2019).
\newblock Convergence complexity analysis of Albert and Chib's algorithm for Bayesian probit regression.
\newblock \textit{The Annals of Statistics, 47}(4), 2320--2347.

\bibitem[{Qin and Hobert(2021)}]{qin:hobe:2021aap}
Qin, Q. and Hobert, J.~P. (2021).
\newblock On the limitations of single-step drift and minorization in Markov chain convergence analysis.
\newblock \textit{The Annals of Applied Probability, 31}(4), 1633--1659.

\bibitem[{Qin and Hobert(2022)}]{qin:hobert:2022aihp}
Qin, Q. and Hobert, J.~P. (2022).
\newblock Geometric convergence bounds for Markov chains in Wasserstein distance based on generalized drift and contraction conditions.
\newblock \textit{Annales de l'Institut Henri Poincar\'e, Probabilit\'es et Statistiques, 58}(2), 872--889.



\bibitem[{Rosenthal(1995)}]{rose:1995}
Rosenthal, J. (1995).
\newblock Minorization conditions and convergence rates for Markov chain Monte Carlo.
\newblock \textit{Journal of the American Statistical Association, 90}(430), 558--566.

\bibitem[{Walczuk(2008)}]{walczuk:2008}
Walczuk, A. (2008).
\newblock Central limit theorem for an additive functional of a Markov process, stable in the Wasserstein metric.
\newblock \textit{Annales Universitatis Mariae Curie-Sklodowska, Sectio A, 62}(1), 149--159.

\end{thebibliography}
\end{document}